\documentclass[12pt,twoside]{article}
\usepackage{amsmath}
\usepackage{amsfonts}
\usepackage{amsthm}
\usepackage{amssymb}
\usepackage{graphicx}
\usepackage{multicol}

\begin{document}
\renewcommand{\baselinestretch}{1.5}
\title{{\normalsize
{\bf Triviality of a surface-link with meridian-based free fundamental group}}}
\author{{\footnotesize Akio KAWAUCHI}\\
\date{}
{\footnotesize{\it Osaka City University Advanced Mathematical Institute}}\\
{\footnotesize{\it Sugimoto, Sumiyoshi-ku, Osaka 558-8585, Japan}}\\
{\footnotesize{\it kawauchi@sci.osaka-cu.ac.jp}}}
\maketitle
\vspace{0.25in}
\baselineskip=15pt
\thispagestyle{empty}
\newtheorem{Theorem}{Theorem}[section]
\newtheorem{Conjecture}[Theorem]{Conjecture}
\newtheorem{Lemma}[Theorem]{Lemma}
\newtheorem{Sublemma}[Theorem]{Sublemma}
\newtheorem{Proposition}[Theorem]{Proposition}
\newtheorem{Corollary}[Theorem]{Corollary}
\newtheorem{Claim}[Theorem]{Claim}
\newtheorem{Definition}[Theorem]{Definition}
\newtheorem{Example}[Theorem]{Example}

\begin{abstract}  
It is proved that every  disconnected surface-link with meridian-based free fundamental group is  a trivial (i.e., an unknotted-unlinked)  surface-link.  
This result is a surface-link version of  the author's recent announcement result 
on smooth unknotting of a surface-knot. 
\end{abstract}

\phantom{x}



\section{Introduction} 

A {\it surface-link}  is a closed oriented 
(possibly disconnected) surface $F$  embedded in the 4-space ${\mathbf R}^4$ 
by a smooth  (or a piecewise-linear locally flat) embedding. 
When $F$ is connected, it is also called a {\it surface-knot}. 
When $F$ is an $r$ copies of the 2-sphere $S^2$, it is called an $S^2$-{\it link} 
with $r$ components.  
For our argument here, a surface-link in the 4-space ${\mathbf R}^4$ is considered as a surface-link in the 4-sphere $S^4$ which is 
the one-point compactification ${\mathbf R}^4\cup\{\infty\}=S^4$. 
Two surface-links $F$ and $F'$  are {\it equivalent} by an {\it equivalence} $f$  if $F$ is sent to $F'$  
orientation-preservingly by an orientation-preserving diffeomorphism  $f:S^4\to S^4$.  
A {\it trivial} surface-link is a surface-link $F$ which is the boundary of the union of mutually disjoint 
handlebodies smoothly embedded in $S^4$, where a handlebody is a 3-manifold which is a 
3-ball or a disk sum of some number of solid tori.  
A trivial surface-knot is also called an {\it unknotted } surface-knot. 
A trivial disconnected surface-link is also called an {\it unknotted-unlinked} 
surface-link. 
For every given closed oriented 
(possibly disconnected) surface $F$,  a trivial  $F$-link exists uniquely  up to equivalences (cf. \cite{HoK}).

The {\it exterior} of  a surface-link $F$ is the compact 4-manifold 
$E=\mbox{cl}(S^4\setminus N(F))$ for a tubular neighborhood $N(F)$ of 
$F$ in $S^4$. Let $q_0$ be a fixed base point  in the interior of $E$. 
A surface-link  with {\it meridian-based free fundamental group} is a surface-link $F$ 
in $S^4$ such that the 
fundamental group $\pi_1(E,q_0)$  is a free group with a meridian basis. 

In \cite{K6} with supplement \cite{KS} generalizing the arguments of ribbon 
surface-links in \cite{K2,K3,K5}, 
the author claimed that 
a surface-knot $F$ in ${\mathbf R}^4$  is a trivial surface-knot if and only if 
it has the infinite cyclic fundamental group. See also \cite{FQ,HiK, K11} for a 
TOP-trivial suface-knot. 

In this paper, 
this result is generalized to any surface-link with  meridian-based free fundamental group. 

Note that there are lots of non-trivial surface-links in  $S^4$ with non-meridian 
based free fundamental groups. For example, 
let $\Gamma$ be a connected non-circular spatial graph 
without degree one vertices  in $S^3$ such that the fundamental group 
$\pi_1(S^3\setminus\Gamma, q_0)$ is a free group and the graph $\Gamma$ has a 
non-trivial constituent knot $\ell$ in $S^3$. Let $B_0$ be a 3-ball in $S^3$ which is a regular neighborhood of a maximal tree of $\Gamma$ in $S^3$, and 
$B=\mbox{cl}(S^3\setminus B_0)$ the complementary 3-ball. Then the intersection $t=\Gamma\cap B$ is a disconnected tangle without loop components in $B$. By Artin's spinning construction of the tangle $(B,t)$, we have a disconnected $S^2$-link $L(t)$ in $S^4$ with the nontrivial 
$S^2$-compnenent obtained by  Artin' spinning construction of the knon-trivial knot 
$\ell$ in $S^3$. As it is observed in \cite[p. 204]{K1},  the fundamental group 
$\pi_1(S^4\setminus L(t),q_0)$ is isomorphic to the fundamental group 
$\pi_1(B\setminus t,q_0)\cong \pi_1(S^3\setminus \Gamma,q_0)$ which is a free group 
by van Kampen theorem. 
Thus, the $S^2$-link $L(t)$ in $S^4$ is a non-trivial $S^2$-link with free fundamental group. 

The following unknotting-unlinking result is our main theorem, answering positively the problem \cite[Problem~1.55 (B)]{Kir} for any $S^2$-link, where note that a trivial surface-link is a surface-link  with  meridian-based free fundamental group.

\phantom{x}

\noindent{\bf Theorem~1.1.}  
A surface-link   is a  trivial surface-link if and only if it is a surface-link 
with meridian-based free fundamental group. 

\phantom{x}

A {\it stabilization} of a  surface-link $F$ is a connected sum 
$\bar F= F\# _{k=1}^s  T_k$ of $F$ and a trivial torus-link $T=\cup_{k=1}^s T_k$.  
By granting $T=\emptyset$, we understand that a surface-link $F$ itself is a stabilization of $F$. 
The trivial torus-link $T$ is called  the {\it stabilizer}  with $T_k\, (k=1,2,...,s)$ the {\it stabilizer components} 
on the stabilization $\bar F$.  

A {\it stably trivial} surface-link  is a surface-link $F$ such that a stabilization 
$\bar F$ of $F$  is 
a trivial surface-link.  In \cite[Corollary~1.2]{K6} with supplement \cite{KS}, it is shown that every stably trivial surface-link is a 
trivial surface-link. 
Therefore, the proof of Theorem~1.1 is completed by combining the result \cite[Corollary~1.2]{K6} with \cite{K6} with supplement \cite{KS}
the following lemma, which generalizes the result  of \cite[Theorem~2.10]{HoK} to a surface-link. 

\phantom{x}

\noindent{\bf Lemma~1.2.}  
Every surface-link  with meridian-based free fundamental group  is a stably trivial surface-link.

\section{Proof of  Lemma~1.2.}

Throughout this section,  the proof of Lemma~1.2 will be given. 
Let $F$ be a surface-link in $S^4$ with components $F_i\, (i=1,2,\dots, r)$. 
 Let $N(F)=\cup_{i=1}^r N(F_i)$ be  a tubular neighborhood of 
$F=\cup_{i=1}^r F_i$ in $S^4$ which is a trivial disk bundle over $F$, and 
$E=\mbox{cl}(S^4\setminus N(F))$  the exterior of  a surface-link $F$. 
Then the boundary 
$\partial E=\partial N(F)=\cup_{i=1}^r \partial N(F_i)$ 
 is a trivial circle bundle over $F=\cup_{i=1}^r F_i$. 
Identify   $\partial N(F_i)=F_i\times S^1$  such that 
the composite inclusion 
\[F_i\times 1\to \partial N(F_i) \to \mbox{cl}(S^4\setminus N(F_i)\]
induces the zero-map in the integral first homology, where $S^1$ denotes  the set of complex numbers of norm one. 
Let   $q_i\times 1$ be a point  in  $F_i\times S^1$ for every $i\, (i=1,2,\dots,r)$. 

Let 
\[K=(\cup_{i=1}^r a_i )\bigcup(\cup_{i=1}^r  S_i)\]
be a connected graph in $E$ such that 

\phantom{x}

\noindent{(1)} $a_i$ is an edge embedded in $E$  joining  $q_0$ and $q_i\times 1$ such that 
the interiors of $a_i\, (i=1,2,\dots,r)$ are mutually disjoint, 

\medskip

\noindent{(2)} $S_i=q_i\times S^1\, (i=1,2,\dots,r)$, and 

\medskip

\noindent{(3)} the inclusion $K\to E$ induces an isomorphism 
$\pi_1(K,q_0)\to \pi_1(E,q_0)$ such that the element $t_i=[a_i\cup S_i]\in \pi_1(E,q_0)$  is the $i$th 
meridian generator.

\phantom{x}

By the assumption that $\pi_1(E,q_0)$ is a free group with a meridian basis, 
there is a graph $K$  with properties (1), (2) and (3).  Further, we have the following lemma. 

\phantom{x}

\noindent{\bf Lemma~2.1.} The composite inclusion $F_i\times 1\to\partial N(F_i) \to E$ is  null-homotopic for all $i$. 

\phantom{x}

\noindent{\bf  Proof of Lemma~2.1.} Since $\partial N(F_i)=F_i\times S^1$, 
the fundamental group elements between the factors $F_i\times 1$ and $q_i\times S^1$ are commutive. 
On the other hand, since $\pi_1(E,q_0)$ is a free group, the image of the homomorphism 
$\pi_1(a_i\cup F_i\times 1 ,q_0)\to \pi_1(E,q_0)$ is in the infinite cyclic group $<t_i>$ generated by $t_i$. 
The surface $F_i\times 1$ in $\partial N(F_i)=F_i\times S^1$  is chosen so that 
the inclusion $F_i\times 1\to\mbox{cl}(S^4\setminus N(F_i))$ 
induces the zero-map in the integral first homology. 
This implies that the  inclusion $F_i\times 1\to E$ is  null-homotopic. \qed 

\phantom{x}

Let 
\[K(\partial N)=(\cup_{i=1}^r a_i )\cup(\cup_{i=1}^r  \partial N(F_i))\] 
be a polyhedron in $E$. 
Let $p: K(\partial N)\to K$ be the map defined by the projection $F_i\times S^1\to q_i\times S^1$ 
sending $F_i$ to $q_i$ for all $i$. 

\phantom{x}

\noindent{\bf Lemma~2.2.} The map $p: K(\partial N)\to K$ extends to a map $g:E\to K$. 

\phantom{x}

\noindent{\bf  Proof of Lemma~2.2.} Since $K$ is a $K(\pi,1)$-space, there is a map $f:E\to K$  inducing 
the inverse isomorphism $\pi_1(E,q_0)\to \pi_1(K,p_0)$ of the isomorphism $\pi_1(K,p_0)\to \pi_1(E,q_0)$. 
Then the composite map $f j: K(\partial N)\to K$ of $f$ with the inclusion $j:K(\partial N)\to E$ and 
the map $p: K(\partial N)\to K$ induces the same homomorphism 
\[(f j)_{\#}=p_{\#}: \pi_1(K(\partial N),q_0) \to \pi_1(K,q_0).\] 
Since $K$ is a $K(\pi,1)$-space, the map $f j$ is homotopic to $p$. 
By the homotpy extension property in \cite{Sp}, there is a map $g:E\to K$ homotopic to the map $f$ 
such that  $g j=p$. \qed

\phantom{x}
 
Replacing $g$ by a piecewise smooth approximation keeping the map $p$ fixed,  we can use a transverse regularity 
argument to obtain 
a regular point $q_i\times t_i\in q_i\times S^1$ for each $i\, (i=1,2,\dots,r)$ such that  
the preimage $V_i = g^{-1}(q_i\times t_i)$ is a compact oriented smooth 3-manifold with boundary 
$\partial V_i= p^{-1}(q_i\times t_i)=F_i\times t_i$. By discarding a closed component from $V_i$, 
we assume that $V_i$ is connected. 
Let $\alpha_{ik}\, (k=1,2,\dots, n_i)$ be mutually disjoint proper arcs in $V_i$ such that 
the exterior $V_i^*=\mbox{cl}(V_i\setminus N_i))$ for a regular neighborhood 
$N_i=N(\cup_{k=1}^{n_i} \alpha_{ik})$ of  the arcs $\alpha_{ik}\, (k=1,2,\dots, n_i)$ in $V_i$ 
is a handlebody.  Let $V(C)_i$ be  a compact connected oriented smooth 3-manifold with boundary 
$\partial V(C)_i=F_i$ obtained from $V_i$ by adding a  collar $C_i$ joining  $F_i\times t_i$  with $F_i$ 
in the trivial disk bundle $N(F_i)$ over $F_i$. 
Let $\alpha(C)_{ik}\, (k=1,2,\dots, n_i)$ be mutually disjoint proper arcs in $V(C)_i$ obtained from 
the arcs $\alpha_{ik}\, (k=1,2,\dots, n_i)$ by extending them until they reach  $F_i$ through the collar $C_i$. 
Then the exterior $V(C)_i^*=\mbox{cl}(V(C)_i\setminus N(C)_i))$ for a regular neighborhood 
$N(C)_i=N(\cup_{k=1}^{n_i} \alpha(C)_{ik})$ of  the arcs 
$\alpha(C)_{ik}\, (k=1,2,\dots, n_i)$ in $V(C)_i$ is a handlebody for all $i$. 
Let $\alpha(C)^0_{ik}$ be an arc obtained from a simple arc $\beta_{ik}$ in $F_i$ with 
 $\partial \beta_{ik}=\partial \alpha(C)_{ik}$ by pushing the interior of $\beta_{ik}$ 
into the collar $C_i$. 
Note that every loop in $V_i$ is null-homotopic in $E$ since 
the map  $g$ induces an isomorphism $g_{\#}:\pi_1(E,q_0)\cong \pi_1(K,q_0)$. 
This means that the arc $\alpha(C)_{ik}$ is homotopic to an arc  $\alpha(C)^0_{ik}$ by a homotopy 
in $(S^4\setminus F)\cup F_i$ relative to $F_i$. By an argument of {HoK}, we see that 
the $1$-handle $h_{ik}$ on $F$ thickening the arc $\alpha(C)_{ik}$ is a trivial 1-handle whose 
surgery of $F$ is the connected sum of $F$ and a trivial torus-knot on the component $F_i$. 
Since the  surface-link obtained from $F$ by the surgery along the 1-handles $h_{ik}$ for all $i$ and $k$  bounds 
a system of handlebodies $V(C)_i^*\, (i=1,2,\dots,r)$, the surface-link $F$ 
is a stably trivial surface-link. 
This completes the proof of Lemma~1.2. \qed

 \phantom{x}

A surface link $F$ with $r(\geq 2)$ components is a {\it boundary surface-link} if it bounds 
mutually disjoint $r$ bounded connected oriented smooth 3-manifolds in $S^4$.  The following corollary 
is contained in the proof of Lemma~1.2 whose technique is  known in classical link theory (see \cite{Sm}). 

\phantom{x}

\noindent{\bf Corollary~2.3.} A disconnected surface link $F$ with $r$ compenents is a boundary surface-link if there is an epimorphism from the fundamental group $\pi_1(E,q_0)$ onto a free group of rank $r$ sending a meridian system to a basis. 

\phantom{x}

\noindent{\bf Acknowledgements.} This work was partly supported by 
Osaka City University Advanced
Mathematical Institute (MEXT Joint Usage/Research Center on Mathematics
and Theoretical Physics JPMXP0619217849). 

\phantom{x}

\end{document}